\documentclass[11pt]{article}
\usepackage{amsmath}
\usepackage{amssymb}
\usepackage{mathptmx}
\usepackage{amsfonts}
\usepackage[mathscr]{eucal}
\usepackage[dvips]{graphicx}
\usepackage{color}

\newcommand{\mbf}[1]{\ensuremath{\mathbf{#1}}}
\newcommand{\mbb}[1]{\ensuremath{\mathbb{#1}}}
\newcommand{\ms}[1]{\ensuremath{\mathscr{#1}}}

\newcommand{\sqzb}[1]{\ensuremath{
d\Gamma_{\textrm{b}}({#1}) }}

\newcommand{\Hb}{H_{\textrm{b}}}
\newcommand{\Hpar}{H_{\textrm{p}}}

\newcommand{\tens}{\otimes}

\newcommand{\nstens}{\otimes^{n}_{\textrm{s}}}
\newcommand{\LtwoRthree}{L^{2}(\mathbb{R}^{3})}

\newcommand{\LtwoRthreeNx}{L^{2}(\mathbb{R}^{3N}_{\,\mbf{x}})}
\newcommand{\LtwoRthreek}{L^{2}(\mathbb{R}^{3}_{ \mbf{k}})}

\newcommand{\Rd}{\mathbb{R}^{d} }
\newcommand{\Rthree}{\mathbf{R}^{3} }
\newcommand{\RthreeN}{\mathbb{R}^{3 N} }
\newcommand{\RthreeNx}{\mathbb{R}^{3 N}_{\,\mbf{x}} }

\newcommand{\Fb}{\mathscr{F}_{\textrm{b}}   } 
\newcommand{\Fbfin}{\mathscr{F}_{\textrm{b}, \textrm{fin}}}

\newcommand{\Omegab}{ \Omega_{\textrm{b}} }
\newcommand{\Pb}{ P_{\, \textrm{b} }}
\newcommand{\PB}{ P_{ \, B} }

\newcommand{\bos}{\textrm{b}}

\newcommand{\phattwo}[2]{\ensuremath{\hat{p}^{#2}_{#1}}}

\newcommand{\I}{\textrm{I}}
\newcommand{\II}{\textrm{II}}

\newcommand{\CILambda}{C_{\textrm{I}}(\Lambda )}
\newcommand{\CIILambda}{C_{\textrm{II}}(\Lambda )}

\newcommand{\Lambdaa}{\Lambda^\alpha}

\newcommand{\rhoLambdaa}{\rho_{\Lambda^\alpha} }

\newcommand{\omegam}{ \omega_{ \, m}}

\newcommand{\rhomLambdaa}{\rho_{m , \Lambda^\alpha} }

\newcommand{\Hbm}{H_{\textrm{b} , \, m}}
\newcommand{\Hkm}{H_{\kappa, m}}
\newcommand{\HkmLambda}{H_{\kappa, m}( \Lambda )}
\newcommand{\HkmLambdaren}{H_{\kappa, m}^{ \textrm{ren}}( \Lambda )}

\newcommand{\HImLambda}{H_{\textrm{I} , \, m}( \Lambda )}
\newcommand{\HImLambdaalpha}{H_{\textrm{I} , \, m}( \Lambda^{\alpha} )}
\newcommand{\corem}{\mathscr{D}_{0, \, m}}

\newcommand{\UkmLambda}{U_{\kappa , m} \left(  \Lambda \right)}
\newcommand{\deltaHkmLambda}{\delta H_{\kappa , m} ( \Lambda)}

\newcommand{\VeffkmLambda}{V_{ \kappa,  \,m}^{\textrm{eff}} ( \Lambda )}
\newcommand{\Veffkm}{V_{\kappa,  m}^{\textrm{eff}} }

\newcommand{\VeffkzeroLambda}{V_{\kappa , 0}^{\textrm{eff}} (\Lambda )}
\newcommand{\Veffkzero}{V_{\kappa , 0}^{\textrm{eff}}}

\newcommand{\VeffkLambda}{V_{\kappa }^{\textrm{eff}} (\Lambda )}
\newcommand{\Veffk}{V_{\kappa }^{\textrm{eff}}}

\newcommand{\1}{{\small \text{1}}\hspace{-0.32em}1}

\newcommand{\bIep}{b_{{}_{\I}}(\epsilon )}
\newcommand{\aII}{a_{{}_{\II}}}
\newcommand{\bII}{b_{{}_{\II}}}

\newtheorem{theorem}{Theorem}[section]
\newtheorem{proposition}[theorem]{Proposition}
\newtheorem{lemma}[theorem]{Lemma}

\newtheorem{remark}{Remark}[section]

\begin{document}
\begin{center}
{\LARGE {Scaling limits   with a removal of ultraviolet cutoffs for    semi-relativistic particles  system coupled to    a scalar Bose field} }\\
$\;$ \\
{ Toshimitsu Takaesu  \\
 \textit{Faculty of Science and Technology, Gunma University,\\ Gunma, 371-8510, Japan } \\
}

\end{center}

\begin{quote}
\textbf{Abstract}.
The system of semi-relativistic particles coupled to a scalar Bose field is considered. 
A scaled  total Hamiltonian for the system  is a self-adjoint operator on a tensor product of a square-integrable function space and a boson Fock space. We consider the strong resolvent limit of a renormalized Hamiltonian, which  is defined  by subtracting a divergent term from the scaled  total Hamiltonian. Applying an abstract scaling limit theory and a unitary transformation, we derive the Yukawa potential and Coulomb potential as  effective potentials. \\

{\small
key words :  Quantum field theory, Relativistic Schr\"odinger operator, Hilbert space.}
\end{quote}

\section{Introduction} 
In this paper,   we consider  a system of  semi-relativistic particles interacting with a scalar Bose field. 
The Hilbert  space for the system is defined  by $
 \ms{H} \, = \, \LtwoRthreeNx \tens \ms{F}_{\bos} (\LtwoRthreek ) $, $N \geq 2$, 
where
$ \ms{F}_{\bos} (\LtwoRthreek )  $ denotes the boson Fock space over $\LtwoRthreek$.   The total Hamiltonian  with a fixed ultraviolet cutoff  $\Lambda_0 > 0$  is defined by
\[ \qquad \quad  
\Hkm  \;  =  \; 
 \Hpar \tens \1 \,
 + \, \1 \tens \Hbm  \, + \,  \kappa H_{\I , \, m} (\Lambda_0 )  , \qquad \qquad \kappa \in \mbf{R} .
\]
Here $\Hpar \,  = \,  \sum\limits_{j=1}^{N}\left( \sqrt{-\triangle_{j} + M^{\,2} \frac{}{} } -M\right) $, $M > 0 $,  is the relativistic Schr\"{o}dinger  operator,    
 $ \Hbm =  \sqzb{\omega_{m}} $   the second quantization of $\omegam (\mbf{k}) = \sqrt{ \mbf{k}^2 + m^2 }$, $m \geq  0$,   the interaction   $
\HImLambda  \, = \,  \sum\limits_{j=1}^N \phi \left( \frac{}{} \frac{\tau_{\, \mbf{x}_j} \rho_{m, \Lambda}}{\sqrt{\omegam}}  \right) 
$  the sum of the scalar field operators where $\tau_{\, \mbf{x}} \rho_{m, \Lambda} (\mbf{k}) \; = \;  \rho_{m, \Lambda} (\mbf{k}) e^{- i \mbf{k} \cdot \mbf{x}} \;   $ and  $ \rho_{m, \Lambda} (\mbf{k})$ is  a characteristic function, and $\kappa \in \mbf{R}$ a  coupling constant.  Let
\begin{equation} \qquad \qquad  \qquad 
\HkmLambda  \;  = \;   \Hpar \tens \1 \,   + \,   \Lambda^2 \,  \1 \tens \Hbm  \, + \,   \kappa \Lambda \, \HImLambdaalpha   ,
 \qquad \qquad  \alpha > 0  ,  \label{5/6.aa}
\end{equation}
where $\Lambda >0$ is a scale parameter. We consider  the asymptotic behaviour of $\HkmLambda $ as $\Lambda \to \infty $ in the strong resolvent sense. The scaling limit for  quantum field Hamiltonians of the type (\ref{5/6.aa}) was  originally  introduced by  Davies \cite{Da79}. The scaled total Hamiltonian is  interpreted as follows. The time evolution of the system  is described  by the one-parameter unitary group $ \{ \exp \left(  - it \HkmLambda \right)\}_{t \in \mbf{R}}$. It  is seen that 
\begin{equation}
\exp \left(  - it \HkmLambda \right)  \, = \,
 \exp \left(  - it_{\Lambda} \left(  \frac{1}{\Lambda^2} \Hpar  \tens \1 + \1 \tens \Hbm + \kappa_{\Lambda} \HImLambdaalpha  \right)\right) ,
\notag 
\end{equation}
where $t_{\Lambda} = t \Lambda^2 $ and $\kappa_{\Lambda} = \frac{\kappa}{\Lambda} $. Formally $ \frac{1}{\Lambda^2 } \left( \sqrt{\mbf{p}^2 + M^{\,2}} -M \right) \simeq 
\frac{1}{2 M_{\Lambda}  } \mbf{p}^2$ for sufficiently small values of momentum $|\mbf{p} | \ll 1$ with
  $M_{\Lambda}= M \Lambda^2 $. 
Thus the investigation of the scaling limit of $\HkmLambda$  as $\Lambda \to \infty$ is    considering  the asymptotic  behaviour of the system with  sufficiently  small  value of the coupling constant   and sufficiently large values of   the time and mass.  
 The scaling limit of $\HkmLambda $ with a fixed ultraviolet cutoff, i.e., $\Lambdaa =\Lambda_{0}$, was  considered  in \cite{Ta11.2}. In this paper, we take  the scaling limit  and remove the  ultraviolet cutoff for the  Hamiltonian simultaneously. Let  $\HkmLambdaren =   \HkmLambda \,  - \, N  E_{\kappa, m} {(  \Lambda)}$ where  $  E_{\kappa , m} {(  \Lambda)}$ is a divergent term as $\Lambda \to \infty $.  Assume  $0 < \alpha < \frac{1}{2}$  and $ 0 <  |\kappa| < 
 \sqrt{\frac{8  \pi}{ N (N-1) }}  $. Then for all $ z \in \mbf{C} \backslash \mbf{R}$,
  \begin{equation}
\textrm{s}-\lim_{\Lambda \to \infty} \, \left(  \HkmLambdaren
  -z \right)^{-1}  \, 
  =  \,   \left(  \Hpar
  \, +  \,   \Veffkm \, - \, z \right)^{-1}  \; \tens \Pb , 
\end{equation}
where  
\[
\Veffkm \; = \;  
-  \,   \frac{\kappa^2}{16\pi}
\sum_{j \, < \,  l}
\frac{1}{|\mbf{x}_{j} -\mbf{x}_{l}|} e^{- m \; |\mbf{x}_{j} -\mbf{x}_{l}| } ,
\]
and $\Pb$ is the projection onto the closed subspace spanned by
 the Fock vacuum $\Omegab$. The potential $\Veffkm$ is called effective potential. We see that    the effective potential  is the Yukawa potential  if $m>0$ and the Coulomb potential if $m=0$.

The scaling limit  with a removal of the ultraviolet cutoff for  a system of  non-relativistic particles  coupled to a scalar Bose fields  was investigated by Hiroshima \cite{Hi97, Hi99} and Gubinelli-Hiroshima-L\"{o}rinczi \cite{GHL14}. The main result in this paper is regarded as   those of   the relativistic version.

 The outline of the proof is as follows. First we investigate  an extension of  the  abstract  scaling limit theory for self-adjoint operators by Arai  \cite{Ar90}. Second we consider  a unitary transformation,  called  \textit{dressing transformation}. 
 From this transformation,  the divergent term, approximating effective potential term and error term are derived. By the abstract scaling limit theory and   the  unitary transformation,   effective potentials are derived.

The scaling limit  for other interacting quantum field models has been analyzed, and refer to    \cite{Hi93, Hi97, Hi02, Oh09, Su07-1, Su07-2, Ta11.1} and reference therein.  The technique of the scaling limit can be applied to  enhanced binding, which is a  phenomenon of ground states appearing. For the enhanced binding of the semi-relativistic particles system coupled to a scalar Bose field, see  Hiroshima-Sasaki \cite{HiSa15}.

This paper is organized as follows. In Section 2,  we consider an abstract
 scaling limit theory for self-adjoint operators.  In Section 3,  we state  the main theorem. In section 4, we prove  the main theorem.


\section{Abstract Scaling Limit } \label{AbstractScalingLimit}
Let  $\ms{Z} = \ms{X} \tens \ms{Y}$ where  $\ms{X}$ and $\ms{Y}$ are  Hilbert spaces. 
Let
\[
X (\Lambda ) \,  =  \, X_0 (\Lambda )   +  \CILambda  +  \CIILambda \tens \1 ,
 \quad \quad \Lambda > 0, 
\]
where  $ X_0 (\Lambda ) \, = \, A \tens \1  +  \Lambda \, \1 \tens B$. Here $A  $ and $\, B  $ are linear operators on $\ms{X}$ and $\ms{Y}$,  respectively,   $\CILambda$ is a linear operator on $\ms{Z}$ and $\CIILambda$ is  a linear  operator on $\ms{X}$. We set 
 $\ms{D}_{A,B} =  \ms{D} (A \tens \1) \cap  \ms{D} (\1 \tens B)  $.
Assume the following conditions.
\begin{quote}
\textbf{[A.1]} $A$  and $B$ are   non-negative and self-adjoint,   and   ker$B \; \ne \{ 0\}$. \\ 
$\;$ \\
\textbf{[A.2]}\textbf{ (i)} $ C_{\I} (\Lambda )$ is a symmetric operator on $\ms{Z}$. For all $\epsilon >0 $ there exists a constant 
 $\Lambda_{\I} (\epsilon ) >0 $ such that for all 
$\Lambda > \Lambda_{\I} (\epsilon )$, it holds that 
 $  \ms{D}_{A,B}  \subset \ms{D} (C_{\I}(\Lambda ))$ and
 \begin{equation} 
  \| \CILambda \Psi  \| \,  \leq \,  \epsilon \| X_0 (\Lambda ) \Psi \| \; + \bIep  \| \Psi \|  ,  \qquad  \Psi \in 
\ms{D}_{A,B} . \notag
 \end{equation} 
where $ \bIep \geq 0 $ is a  constant independent of $\Lambda > \Lambda_{\I} (\epsilon )$. \\
 \textbf{(ii)} Moreover,  
there exists a symmetric operator $C_{\I }$ on $\ms{Z}$ such that
$ \ms{D}_{A,B}  \subset \ms{D} (C_{\I}) $, and 
\begin{equation}
\text{s}-\lim_{\Lambda \to \infty } \CILambda  (X_0 (\Lambda ) -z)^{-1}
  \; = \;  C_{\I } \left( (A-z)^{-1}  \tens \PB  \right) , \qquad z \in  \mbb{C} \backslash \mbb{R}  ,   \notag
\end{equation}
where $\PB $ is the projection onto ker$B$. \\
$\;$ \\
\textbf{[A.3]} \textbf{(i)} $ C_{\II} (\Lambda )$ is a symmetric operator on \ms{X}. For all $\Lambda > 0$, $\ms{D} (A)  \subset \ms{D} (\CIILambda )$ and there exist constants $  0 <  \aII <1  $ and $ \bII \geq 0$ such that for all  $\Lambda > \Lambda_{\II}$ with some $\Lambda_{\II} >0 $,
 \begin{equation}
  \| \CIILambda \psi  \| \; \leq \;
 \aII \| A \psi \| \; +  \bII  \| \psi \| ,  \qquad \psi \in \ms{D} (A) .  \notag
 \end{equation}
\textbf{(ii)} In addition, there exists a symmetric operator $C_{\II  }$ on $\ms{X}$ such that $\ms{D} (A) \subset \ms{D} (C_{\II}) $,
\begin{equation}
\text{s}-\lim_{\Lambda \to \infty } (\CIILambda  \tens \1) (X_0 (\Lambda ) -z)^{-1}
 \;  = \;   \left( C_{\II } (A-z)^{-1}  \right) \tens \PB , \qquad z \in  \mbb{C} \backslash \mbb{R}   . \notag
\end{equation}
 \end{quote}

\begin{proposition} \upshape \label{AbstractScaling}
Assume \textbf{[A.1]} - \textbf{[A.3]}. Let $\Lambda_{0}(\epsilon) = \textrm{max} \{ \Lambda_{\I}(\epsilon) ,\Lambda_{\II} \}$. Then,   \\
\textbf{(i)}  Let $\epsilon >0 $ such that  $\aII + \epsilon <1$. Then for all $\Lambda > \Lambda_{\, 0}(\epsilon)$,
$X (\Lambda )$ is self-adjoint on $\ms{D}_{A,B}$ and essentially self-adjoint on any core of $X_0 (\Lambda) $. \\
\textbf{(ii)} The operator 
\[
X_\infty   \; = \;  A\tens \1  \, + \,  ( \1 \tens \PB ) \, 
  C_{\I}   ( \1 \tens \PB)  \, + \,   C_{\II} \tens  \PB 
\]
is self-adjoint on $\ms{D} (A \tens \1)$ and essentially self-adjoint on any core of $A \tens \1 $.\\
\textbf{(iii)}
For all $z \in  \mbb{C} \backslash \mbb{R} $,
\[
\text{s} - \lim_{\Lambda \to \infty } \left(  X (\Lambda ) -z \right)^{-1}
\;  = \; ( X_{\infty}  \; - \; z  )^{-1} ( \1 \tens \PB  ). 
\]
\end{proposition}

\begin{remark} \upshape $\; $ 
Proposition \ref{AbstractScaling}  with  $C_{\II} (\Lambda )=0$  was proven by  Arai  (\cite{Ar90}; Theorem 2.1).
\end{remark}

$\quad $ \\
\textbf{{\large (Proof)}}  
\textbf{(i)} Let $\Psi \in   \ms{D}_{A,B}$. 
   It holds that $ \| ( A \tens \1 ) \Psi \|  \leq \left\| X_{0} (\Lambda ) \Psi  \right\|$. By \textbf{[A.2]} and \textbf{[A.3]} it follows  that for  all $\Lambda > \Lambda_{\, 0}(\epsilon)$,  
\begin{equation}
\left\| \left(  \CILambda + \CIILambda \tens \1 \right) \Psi \right\|
 \leq ( \aII + \epsilon ) \|   X_{0} (\Lambda ) \Psi  \|
  + (\bII + \bIep   ) \| \Psi \| . \label{5/8.7}
\end{equation}
Note that    $X_0 (\Lambda ) $ is self-adjoint on $\ms{D}_{A ,B}$, since $A$ and $B$ are  self-adjoint and bounded from below. Taking    $\epsilon > 0$ such as  $ \aII + \epsilon <1$, we  obtain  \textbf{(i)}  from  the Kato-Rellich theorem. \\
$\quad$ \\
\textbf{(ii)} Let $\Xi \in   \ms{Z}$. 
By (\ref{5/8.7}),  we have
\begin{align}
&  \left\| \left(  \CILambda + \CIILambda \tens \1  \right)  \left( \frac{}{} X_0 (\Lambda ) -z \right)^{-1} \Xi \right\|  \notag \\
&  \qquad   \leq \,  ( \aII + \epsilon ) \left\| X_{0} (\Lambda ) \left( X_{0} (\Lambda )  -z \right)^{-1} \Xi  \right\|    +  \left(  \bII + \bIep \right) 
\| \left(X_{0} (\Lambda )  -z \right)^{-1} \Xi  \|  \notag  \\
& \qquad \leq \,  ( \aII + \epsilon ) \| \Xi \| + \left( ( \aII+  \epsilon) |z| +\bII +  \bIep   \right) 
\left\| \left(X_{0} (\Lambda )  -z \right)^{-1} \Xi  \right\|  . \notag 
\end{align}
Note that s-$\lim\limits_{\Lambda \to \infty}  \left(X_{0} (\Lambda )  -z \right)^{-1} = (A-z)^{-1} \tens \PB  $.  From the above inequality,   we have
\begin{align}
&\left\| \left(  C_{\I} + C_{\II} \tens \1    \right)  \left( \frac{}{} (A-z)^{-1}  \tens \PB \right) \Xi \right\|   \notag \\
& \qquad \qquad \leq  ( \aII +\epsilon ) \|  \Xi  \| + \left( \frac{}{} ( \aII + \epsilon ) |z| + \bII + \bIep  \right) 
\left\|  \left( \frac{}{} (A-z)^{-1}  \tens \PB \right) \Xi  \right\| .  \notag
\end{align}
Let $ \Xi = ( (A-z) \tens \1 )\Psi $, $\Psi \in \ms{D} (A \tens \1 )$. Then  we have  
\begin{align}
&\left\|  ( \1 \tens \PB ) \left( C_{\I} + C_{\II} \tens \1  \right)  ( \1 \tens \PB ) \Psi \right\| \notag \\
& \qquad \qquad   \leq   ( \aII + \epsilon ) \| ( A\tens \1   ) \Psi \| +\left( \frac{}{}  2 ( \aII + \epsilon ) |z|  + \bII + \bIep  \right) \|  \Psi \| . \notag
\end{align}
Taking  $\epsilon >0 $  such that $ \aII + \epsilon  < 1$, we obtain \textbf{(ii)}  from the Kato-Rellich theorem. 

$\; $ \\
\textbf{(iii)}
 Let   $z \in  \mbb{C} \backslash \mbb{R} $ with Re $z$ $<0$. 
  From (\ref{5/8.7}), it follows that  for all $\Lambda > \Lambda_{\, 0} (\epsilon )$,
 \begin{equation}
\left\| \left( \CILambda + \CIILambda \tens  \1  \right)  \left( X_0 (\Lambda ) -z \right)^{-1} \right\| 
 \leq  ( \aII + \epsilon  ) + \frac{\bII + \bIep  }{|z|} . \notag
\end{equation}
Let $\epsilon >0 $  such that $ \aII + \epsilon  < 1$ and assume that $|z|>0$ is  sufficiently large.  Then it holds that $\| (  \CILambda + \CIILambda \tens \1 )  (  X_0 (\Lambda ) -z )^{-1} \| 
 \, < \,  1$.  We see that s-$\lim\limits_{\Lambda \to \infty} ( X_0 (\Lambda) - z )^{-1} = (A-z)^{-1} \tens \PB$. Then  by the Neumann series expansion, we have 
\begin{align*}
\text{s-}\lim_{\Lambda \to \infty }\left( X ( \Lambda ) -z \right)^{-1} 
& = \text{s-}\lim_{\Lambda \to \infty }  \sum_{n=0}^\infty (-1)^n \left( X_0 ( \Lambda ) -z \right)^{-1} 
\left\{   \left( \CILambda + \CIILambda \tens \1  \right) 
 \left(  X_0 (\Lambda ) -z \right)^{-1}
 \right\}^n  \\
 &= \sum_{n=0}^\infty  (-1)^n \left( \frac{}{} ( A -z )^{-1} \tens \PB \right)
\left\{ \frac{}{}  \left(  C_{\I} + C_{\II} \tens \1  \right) \left(  \left( A -z \right)^{-1}  \tens \PB \right)
 \right\}^n  \\
 & =(X_{\infty }-z )^{-1} ( \1 \tens \PB ) .
\end{align*}
Thus  \textbf{(iii)} is proven. $\blacksquare$ \\

To prove the main theorem, it is useful to prepare the following lemma.

\begin{lemma} \label{A.3'} \upshape
Assume \textbf{[A.1]} and  \textbf{[A.3]} \textbf{(i)}. In addtion, suppose that   
there exists a symmetric operator $C_{\II  }$ on $\ms{X}$ such that $\ms{D} (A) \subset \ms{D} (C_{\II}) $ and for all $\psi \in \ms{D} (A)$,
\begin{equation}
\text{s}-\lim_{\Lambda \to \infty } \CIILambda  \psi
 \;  = \;   C_{\II} \psi  . \notag 
\end{equation}
Then \textbf{[A.3]} \textbf{(ii)} holds. 
 \end{lemma}
\textbf{(Proof)}  Let $z \in \mbb{C} \backslash \mbb{R}$. 
 We see that 
\[
 ( \CIILambda \tens \1 ) ( X_0 (\Lambda) - z )^{-1} = \left( \CIILambda (A-z)^{-1} \right) \tens \PB +  ( \CIILambda \tens \1 )( X_0 (\Lambda) - z )^{-1} (\1 \tens \PB^{\bot}) .
\]
Then it holds that for all $\Psi \in \ms{Z}$,
\begin{align}
& \lim_{\Lambda \to \infty} \| ( \CIILambda  \tens \1 )\left(  X_0 (\Lambda ) -z \right)^{-1} (\1 \tens \PB^{\bot} ) \Psi  \|  \notag \\
 & \; \; \leq    \lim_{\Lambda \to \infty} \left( 
\aII \|  ( A  \tens \1 ) (X_0 (\Lambda ) -z )^{-1} ( \1 \tens \PB^{\bot} ) \Psi \|  + \bII   \| \left(  X_0 (\Lambda ) -z \right)^{-1} ( \1 \tens \PB^{\bot} ) \Psi \|  \right) \notag \\
& \;\; = \lim_{\Lambda \to \infty} \left( 
\aII \| \left( X_0 (\Lambda ) -z \right)^{-1} ( \1 \tens \PB^{\bot} ) ( A  \tens \1 ) \Psi  \|  + \bII  \| \left(  X_0 (\Lambda ) -z \right)^{-1} ( \1 \tens \PB^{\bot} ) \Psi \|  \right)
= 0  \notag .
\end{align}
From this,    we have  s-$\lim\limits_{\Lambda \to \infty } ( \CIILambda \tens \1 )\left(  X_0 (\Lambda ) -z \right)^{-1} \Psi
 \,  = \,   \left( C_{\II} (A-z)^{-1} \tens \PB \right) \Psi $. $\blacksquare$

\section{Main Theorem}
The definition of the  state space and total Hamiltonian for the semi-relativistic particles interacting with a scalar Bose field is as follows. The state space is defined by  
\[
\qquad \qquad \qquad 
 \ms{H} \; = \; \LtwoRthreeNx \tens \ms{F}_{\bos} (\LtwoRthreek ) , \qquad \qquad \qquad N \geqq 2. 
\]
where
$ \Fb (\LtwoRthreek ) \, = \,  \oplus_{n=0}^{\infty}
 (  \nstens \LtwoRthreek )  $ is the boson Fock space over $\LtwoRthreek$. Here  $  \nstens  \ms{X}$ denotes the $n$-fold symmetric tenser  product of a Hilbert space $\ms{X}  $ and $\tens_{s}^{0} \ms{X}= \mbb{C}$. 
The total Hamiltonian  with the fixed ultraviolet cutoff  $\Lambda_0 > 0$  is defined by
\[ \qquad \quad  
\Hkm  \;  = \; 
 \Hpar \tens \1 \,
 + \, \1 \tens \Hbm  \, + \,   \kappa \, H_{\I , \, m} (\Lambda_0 )  , \qquad \qquad \kappa \in \mbf{R}.
\]
The  semi-relativistic particle Hamiltonian is defined by $\Hpar \,  = \,  \sum\limits_{j=1}^{N}\left( \sqrt{ -\triangle_{j} + M^{\, 2} \frac{}{}  } -M\right)  $, $M>0$, and 
the  field Hamiltonian by 
 $ \Hbm =  \sqzb{\omegam } $ with $\omegam (\mbf{k}) = \sqrt{  \mbf{k}^2 + m^{\, 2} }$, $m \geqq 0$, where $\sqzb{X }$ denotes the second quantization of the operator $X$. The interaction  is defined by
\begin{equation}
\HImLambda  \; = \;  \sum_{j=1}^N \phi \left( \frac{}{} \frac{\tau_{\, \mbf{x}_{j}}\rho_{m,\Lambda }}{\sqrt{\omegam}}  \right) , \notag 
\end{equation}
where $(\tau_{\, \mbf{x}} f) (\mbf{k}) \, = \,  f 
(\mbf{k}) e^{- i \mbf{k} \cdot \mbf{x}} \;   $, $f \in \LtwoRthree$, and   $\rho_{m, \Lambda} 
(\mbf{k}) = \chi_{{}_{J_{m, \Lambda}}} (|\mbf{k}|)$ is a characteristic function on 
\[
J_{m, \Lambda} \; = \; 
\left\{
\begin{array}{ll}
[ 0, \,  \Lambda ] &  \quad m>0 , \\
 {[} \frac{1}{\Lambda}, \, \Lambda  {]}  & \quad   m=0 .
\end{array}
\right.
\]
Here $\Lambda >0 $ denotes the ultraviolet cutoff and $\Lambda^{-1}$ the infrared cutoff.
The field operator $\phi (f)$ is defined by $\phi (f)   =   \frac{1}{\sqrt{2}} \left( \frac{}{} a (f )  + a^{\ast} (f) \right) $, $f \in \LtwoRthreek $, where
       $a^{\ast}(f)$ and 
 $a(f)$ denote the creation operator and the annihilation operator, respectively.  The creation operators and the  annihilation operators  satisfy the canonical commutation relations : 
\begin{equation}
 [a (f) , a^{\ast} (g )] = (f ,g ), \qquad \quad 
 [a(f) , a (g)] = [a^{\ast} (f) , a^{\ast} (g)] = 0  , \label{CCR}
\end{equation}
 on the finite particle 
subspace $\Fbfin (\LtwoRthreek)$. 
 The finite particle subspace 
 $\Fbfin (\ms{M})$ on the subspace $\ms{M} \subset \LtwoRthreek $ is defined by the set  which consists of all  vectors of the form  $\Psi = \{\Psi^{(n)} \}_{n=0}^{\infty} \; $ such  that
   $\Psi^{(n)}  \in \nstens \ms{M}  $, $ n  \geq 0 $, and
    $\Psi^{(n')}  = 0  $ for all $ n' > n_0 $ with  some $n_0 \geq 0$.
 
 The self-adjointness of $\Hkm $ is as follows. 
  Let $\omega$ be a non-negative function on $\Rd $, $d \geq 1$, and   $f  \in \ms{D} (\omega^{-1/2})  $. Then for all $\Psi \in  \ms{D}(\sqzb{\omega}^{1/2})$,  
\begin{align} & \| a (f) \Psi  \| \leq \| \frac{f}{\sqrt{\omega}}  \|  \, \|  \sqzb{\omega}^{1/2} \Psi \| ,        \label{bounda}   \\  & \| a^{\ast}(f ) \Psi  \| \leq \|  \frac{f}{\sqrt{\omega}} \| \, \|   \sqzb{\omega}^{1/2}  \Psi \| + \| f \| \| \Psi \|  . \label{boundad} 
\end{align}
From (\ref{bounda}) and (\ref{boundad}), it is  shown that  $H_{\I , \, m} (\Lambda_0 )$ is relatively bounded with respect to $\1 \tens \Hbm^{1/2}$. Then it is easy to show that  $H_{\I , \, m} (\Lambda_0 ) $ is relatively bounded with respect to $H_{0 , \, m}$ with infinitely small bound, and hence,  $\Hkm$ is self-adjoint and essentially self-adjoint on any core of $H_{0 , \, m}$ from the Kato-Rellich theorem. 
 In particular, $\Hkm $ is essentially self-adjoint on 
\[
  \corem \;  =  \; C_{0}^{\,\infty} (\RthreeNx )  \hat{\tens} \Fbfin (\ms{D} (\omegam )) ,
\]
 where $\hat{\tens}$ denotes the algebraic tensor product.

We define a scaled total Hamiltonian by 
\begin{equation} \qquad \qquad  \qquad 
\HkmLambda  \, = \,  H_{0 , \, m } (\Lambda )  +   \kappa \Lambda \, \HImLambdaalpha   ,
 \qquad \qquad  \alpha > 0  ,  \notag
\end{equation}
where $H_{0 , \, m} (\Lambda ) \; = \; \Hpar \tens \1    +    \Lambda^2 \, ( \1 \tens \Hbm ) $.
A unitary transformation, called  \textit{dressing transformation}, is defined by 
\[
  \UkmLambda \;  =  \; \exp \left( {i \left( \frac{\kappa}{\Lambda} \right) \, \sum\limits_{j=1}^N  \pi 
 \left( \frac{ \tau_{\, \mbf{x}_j} \rho_{m, \Lambda^{\alpha}  }}{ \sqrt{\omegam}^3 }  \right) } \right) , 
\]
where $  \pi (f) \;  = \;  \frac{i }{\sqrt{2}} \left( \frac{}{} -  a (f )  + a^{\ast} (f) \right) $, $f \in \LtwoRthreek $, is the conjugate operator. It holds that
\begin{align*}
& [ \pi (f) , \Hbm ] \; = \; -i \,   \phi ( \omegam  f )  ,   \qquad  f \in \ms{D} (\omegam ), \\
& [  \pi (g) , \phi (h ) ]\; =  \; \,  - i \text{Re} (g , h ) ,  \qquad   g, h \in \LtwoRthreek ,  
\end{align*}
 on  $\Fbfin (\LtwoRthreek)$. It is seen that
\begin{equation}
\UkmLambda^{-1} 
 \HkmLambda \UkmLambda
   \; = \;  H_{0} (\Lambda )  \, + \,  N  E_{\kappa , m} ( \Lambda ) \,  + \, \deltaHkmLambda  \, + 
     \,  \VeffkmLambda      , \label{u-transformed}
\end{equation}
where 
\begin{align*}
&E_{ \kappa ,m } ( \Lambda) \; \; = \; \; \frac{\kappa^2}{(2\pi )^3}  \int_{\Rthree} \frac{\rho_{m, \,\Lambda^{\alpha}}(\mbf{k})}{
  \omegam ( \mbf{k} )^2 }  d 
\mbf{k} , \\
& \deltaHkmLambda \; = \; 
 \UkmLambda^{-1}
  \left( \Hpar  \tens \1 \right) \UkmLambda \;  -\;  \Hpar \tens \1   ,  \\
  &  \VeffkmLambda \;  = -  \frac{\kappa^2}{2(2\pi )^3} \; \, \sum_{j \, < \,  l }^N \; \int_{\Rthree}  \frac{ 
\rho_{m , \, \Lambda^\alpha }(\mbf{k} )}{\omegam (\mbf{k})^2} 
 e^{i \mbf{k} \cdot (  \mbf{x}_{j} -\mbf{x}_{l}) }
 d \mbf{k} .
\end{align*}
By  the abstract scaling limit theory in Section \ref{AbstractScalingLimit} and the unitary transformed Hamiltonian  (\ref{u-transformed}), we prove the following  theorem.
 
\begin{theorem} \upshape \label{MainTheorem} 
Assume  $0 < \alpha < \frac{1}{2}$  and $ 0 <  |\kappa| < 
 \sqrt{\frac{8  \pi}{ N (N-1) }}  $. Then for all $ z \in \mbf{C} \backslash \mbf{R}$,
  \begin{equation}
\textrm{s}-\lim_{\Lambda \to \infty} \, \left( \HkmLambdaren
  -z \right)^{-1}     \; 
  =  \;  \left(  \Hpar
   +    \Veffkm  -  z \right)^{-1}  \tens \Pb , \notag
\end{equation}
where $\HkmLambdaren =   \HkmLambda \,  - \, N  E_{\kappa ,m} {( \Lambda)}$ and    
\[
\Veffkm (\mbf{x}_1 , \cdots , \mbf{x}_{N} )\, = \, 
-  \frac{\kappa^2}{16\pi}
\sum_{j \, < \,  l}
\frac{1}{|\mbf{x}_{j} -\mbf{x}_{l}|} e^{- m \; |\mbf{x}_{j} -\mbf{x}_{l}| } ,
\]
and $\Pb$ is the projection onto the closed subspace spanned by the Fock vacuum $\Omegab$. 
\end{theorem}

\section{Proof of Theorem \ref{MainTheorem}}
We show that $  H_{0} (\Lambda ) $, $ \deltaHkmLambda $  and $\VeffkmLambda$
  satisfy the   conditions \textbf{[A.1]} - \textbf{[A.3]} in Section 2. 
\begin{proposition} \upshape \label{5/7.5}
 Assume $0 < \alpha  < \frac{1}{2} $. Then the following holds.   \\
 \textbf{(i)} For all $\epsilon >0 $, there exists $\Lambda_{\, \ast} (\epsilon) >0 $ such that for all $\Lambda >  \Lambda_{\, \ast} (\epsilon)  $,
 \begin{equation}
 \qquad \qquad 
  \| \deltaHkmLambda \Psi  \| \; \leq \; 
   \epsilon \left( \| H_{0 , \, m } (\Lambda ) \Psi \|  + \| \Psi \| \frac{}{}\right) ,  \quad \quad  \Psi \in \ms{D} 
(H_{0 , \, m} (\Lambda )) . \notag
 \end{equation}
\textbf{(ii)} It follows that 
\begin{equation}
\quad \quad 
\text{s}-\lim_{\Lambda \to \infty } \deltaHkmLambda ( H_{0 , \, m } (\Lambda ) -z )^{-1} 
\; = \;  0 ,   \qquad \qquad  z \in \mbb{C} \backslash \mbb{R} .  \notag
\end{equation}
\end{proposition}
$\quad$ \\
\textbf{{\large (Proof)}}  
 \textbf{(i)} 
Let $W_1 $ and  $W_2 $  be  non-negative and self-adjoint operators with   $\ms{D}(W_2) \subset  
 \ms{D} (W_1) $. It is seen  that for all $ \psi \in \ms{D}(W_1) \cap \ms{D}(W_2) $,
 $ \sqrt{W_j} \psi \; = \; 
 \frac{1}{\pi} \int_{0}^{\infty}
 \frac{1}{\sqrt{\lambda}} (W_j + \lambda  )^{-1} \,  W_j \,  \psi
 \;  d \lambda $,$j= 1,2$. Then, we have $
 (\sqrt{W_{1}} - \sqrt{W_{2} } ) \psi = \frac{1}{\pi} \int_{0}^{\infty} 
 \sqrt{\lambda} (W_{1} + \lambda )^{-1} (W_1 - W_2) (W_{2} + \lambda )^{-1} \psi d \lambda  $, and  
 \begin{equation}
 \| (\sqrt{W_{1}} - \sqrt{W_{2} } ) \psi  \| \leq  \frac{1}{\pi} \int_{0}^{\infty} \frac{\sqrt{\lambda}}{\lambda + E_{0} (W_1 ) }  \left\| \frac{}{}  (W_1 - W_2) (W_{2} + \lambda )^{-1} \psi \right\| d \lambda  ,  \label{5/8.3}
\end{equation}
 where
 $E_{0} (X) $ is the infimum of the spectrum of the operator $X$. Let  $\hat{\mbf{p}}_j  $ $= 
 ( \hat{p}_j^\nu  )_{\nu =1}^3  = 
 (   -i \frac{\partial}{ \partial x^\nu_j }  )_{\nu=1}^3 $, $j=1 , \cdots,N $.
By the commutation relation  $[ \pi ( \tau_{\, \mbf{x}_{j}} f ) , \, \hat{p}_{l}^{\nu}  ] \; = \; i \pi (\partial_{x_{l}^\nu} (\tau_{\, \mbf{x}_{j}} f) ) \delta_{j,l} $, $f \in \LtwoRthree$, 
 it holds that  for all $ \Phi \in \ms{D}(H_{0,m } (\Lambda)) $, 
\begin{equation}
\left( \UkmLambda^{-1} ( \hat{\mbf{p}}_{j} \tens \1 ) \UkmLambda  \right)^2  \Phi  \;
 =  \; 
   \left( \sum_{\nu=1 }^{3} 
\left( \phattwo{j}{\nu}  \tens \1  + \left( \frac{\kappa}{\Lambda} \right)
 \pi \left( \frac{  \partial_{x_{j}^\nu } ( \tau_{\, \mbf{x}_{j}} \rho_{m, \Lambda^{\alpha}} )}{ \omegam^{3/2}} \right)
\right)^2  \;  \right) \Phi .  \notag 
\end{equation}
Applying    $\UkmLambda^{-1} (\hat{\mbf{p}}_j^2 \tens \1 ) \UkmLambda  + M^{\, 2}$ to $W_1$,  and  $\hat{\mbf{p}}_j^2 \tens \1+M^{\, 2}$  to  $W_2$  in (\ref{5/8.3}), we have
for all  $\Lambda > \Lambda_{\, \star }$  with $\Lambda_{\, \star } \geq 1 $ and for all $ \Psi \in \ms{D}(H_{0,m } (\Lambda)) $,
\begin{align}
&\| \deltaHkmLambda \Psi \|  \notag \\
    & \; \; \leq \sum_{j=1}^N  \left\| \left( \sqrt{\left( \UkmLambda^{-1} ( \hat{\mbf{p}}_{j} \tens \1 ) \UkmLambda  \right)^2  +M^{\, 2}} 
- \sqrt{\left(   \hat{\mbf{p}}_{j}^2 \tens \1 \right)     +M^{\, 2}} \right) \Psi  \right\| \notag \\
& \; \; \leq \frac{1}{\pi} \sum_{j=1}^N \int_{0}^{\infty} \frac{\sqrt{\lambda}}{\lambda + M^{\,2} }
  \left\| \left(  \left(  \UkmLambda^{-1} (\hat{\mbf{p}}_j \tens \1 ) \UkmLambda   \right)^2 
- \hat{\mbf{p}}_j^2 \tens \1 \right) ( \hat{\mbf{p}}_{j}^2 \tens \1  \; + \; M^{\,2}  + \lambda )^{-1} \Psi \right\|
   \notag \\
&\; \; \leq \frac{1}{\pi} \left( \frac{\kappa}{\Lambda} \right) 
 \sum_{j=1}^N \int_{0}^{\infty} \frac{\sqrt{\lambda}}{\lambda + M^{\,2} }
 \left\| \frac{}{} S_{m, j} (\Lambda) ( \hat{\mbf{p}}_{j}^2 \tens \1  \; + \; M^{\,2}  + \lambda )^{-1} \Psi \right\| \;
d \lambda  \notag \\
 &\qquad \qquad  +\frac{1}{\pi} \left( \frac{\kappa}{\Lambda} \right)^2
 \sum_{j=1}^N \int_{0}^{\infty} \frac{\sqrt{\lambda}}{\lambda + M^{\, 2} } \left\| 
 T_{m, j} (\Lambda) (  \hat{\mbf{p}}^2_{j} \tens \1 \; + \; M^{\, 2}  + \lambda )^{-1} \Psi \right\|  \;  d \lambda ,
  \label{5/7.1}
\end{align}
where 
\begin{align*}
 &S_{m, j} (\Lambda) \;\; = \; \;
-2 \sum_{\nu=1}^3 \pi \left( \frac{i k^\nu \tau_{\, \mbf{x}_j} \rhomLambdaa }{ \omegam^{3/2}} \right)
  ( \hat{p}_j^\nu \tens \1 ) \;  - \;  \pi \left( \frac{ \mbf{k}^2 \tau_{\, \mbf{x}_j}\rhomLambdaa }{ \omegam^{3/2}} \right) ,  \\
 &T_{m, j} (\Lambda) \;\; = \; \; \sum_{\nu=1}^3 \pi \left( \frac{i k^\nu \tau_{\, \mbf{x}_j} \rhomLambdaa }{ \omegam^{3/2}} \right)^2 .
\end{align*}
 By the following Lemma \ref{5/7.2}, it holds that {\small
\begin{align}
& \left\|  S_{m, j} (\Lambda) (  \hat{\mbf{p}}^2_{j} \tens \1 \, + \, M^{\, 2} 
 + \lambda )^{-1} \Psi \right\| 
  \; \leq \;  
\left( 2 \beta_j  \frac{\Lambda^\alpha}{(\lambda + M^{\,2})^{\frac{3}{4}}} + \gamma_j \frac{\Lambda^{2 \alpha}}{\lambda + M^{\, 2}}\right) \left( \left\|  H_{0 , \, m} (\Lambda ) \Psi \right\|   + \| \Psi \| \right)   ,  \label{5/2.1} \\
&\left\|   T_{m, j} (\Lambda) (   \hat{\mbf{p}}^2_{j} \tens \1  \, + \, M^{\, 2}  + \lambda )^{-1} \Psi \right\| 
 \;  \leq  \;  \mu_j \frac{\Lambda^{2\alpha}}{\lambda + M^{\, 2}}
 \left( \left\| H_{0 , \, m } (\Lambda ) \Psi \right\|   + \| \Psi \| \right)   ,  \label{5/2.2} 
\end{align}}
where $\beta_j  \geq 0 $,  
 $\gamma_j \geq 0 $ and $\mu_j \geq 0 $ are constants. 
By (\ref{5/7.1}), (\ref{5/2.1}) and (\ref{5/2.2}), we have
\begin{equation}
\qquad 
\| \deltaHkmLambda \Psi \| \; \; 
 \leq  \left( \; \;  \frac{C_1}{\Lambda^{1-\alpha}} \,
 +  \, \frac{C_2}{\Lambda^{1-2\alpha}} \, 
  +  \, \frac{C_3}{\Lambda^{2(1-\alpha )}}   \right)    \left( \left\|  H_{0 , \, m } (\Lambda ) \Psi \right\|   + \| \Psi \| \right)  ,    \notag
\end{equation}
where $C_1 = \sum\limits_{j=1}^N \frac{\kappa \beta_j }{\pi}  
\int_{0}^{\infty} \frac{ \sqrt{\lambda}}{(\lambda + M^2 )^{ \frac{7}{4}} } d \lambda$, 
 $C_2 = \sum\limits_{j=1}^N  \frac{\kappa \gamma_j }{\pi}  
 \int_{0}^{\infty} \frac{\sqrt{\lambda }}{(\lambda + M^2 )^2} d \lambda $, and $C_3 = \sum\limits_{j=1}^N\frac{\kappa^2 \mu_j }{\pi}   \int_{0}^{\infty} \frac{\sqrt{\lambda}}{(\lambda + M^2 )^2} d \lambda $. 
 Since $\Lambda \geq \Lambda_{\, \star}$, we have 
 \[
 \| \deltaHkmLambda \Psi \| \; \; 
 \leq \epsilon (\Lambda_{\, \star})  \left( \left\|  H_{0 , \, m } (\Lambda ) \Psi \right\|   + \| \Psi \| \right) ,
 \]
 where $\epsilon (\Lambda_{\, \star}) =\frac{C_1}{\Lambda_{\, \star}^{1-\alpha}} \,
 +  \, \frac{C_2}{\Lambda_{\, \star}^{1-2\alpha}} \, 
  +  \, \frac{C_3}{\Lambda_{\, \star}^{2(1-\alpha )}} $. 
 By taking sufficiently large $\Lambda_{\, \star}> 0$, we obtain  \textbf{(i)}.  \\
$\;$ \\
\textbf{(ii)} Let  $z \in \mbb{C} \backslash {\mbb{R}}$ and $\epsilon >0$.  From \textbf{(i)}, it is seen that 
for all $\Psi \in \ms{D} ( H_{0 , \, m } (\Lambda )) $ and  $\Lambda >\Lambda_{\, \ast} (\epsilon)  $ ,
 \begin{align}\| \deltaHkmLambda (H_{0 , \, m } (\Lambda ) -z)^{-1}\Psi  \| \; & \leq \; 
   \epsilon \left( \| H_{0 , \, m } (\Lambda )(H_{0 , \, m } (\Lambda ) -z)^{-1} \Psi \|  + \|(H_{0 , \, m } (\Lambda ) -z)^{-1} \Psi \| \frac{}{}\right) \notag \\
& \leq \epsilon \left(1 + \frac{1}{|\text{Im} \,z |}  \right) \| \Psi \| . \notag
\end{align}
Since $\epsilon >0$ is arbitrary, we obtain \textbf{(ii)}. $\blacksquare $ \\

\begin{lemma} \upshape \label{5/7.2}
Let $\Lambda \geq 1$. 
Then there exist constants $\beta_j \geq 0 $,   
 $\gamma_j \geq 0 $ and $\nu_j \geq 0 $, which are independent of $\Lambda$, and satisfy 
\begin{align*}
&\textbf{(i)} \quad \sum_{\nu =1}^3 \left\|  \pi \left( \frac{ k^\nu \tau_{\, \mbf{x}_j}\rhomLambdaa }{ \omegam^{3/2}} \right)
  ( \hat{p}_j^\nu \tens \1 ) 
 ( \hat{\mbf{p}}_j^2   \tens \1 \; + \; M^{\, 2}  + \lambda )^{-1} (H_{0 , \, m } 
(\Lambda ) +1)^{-1} \right\| 
 \;  \leq  \; 
 \beta_j  \;  \frac{\Lambda^\alpha}{(\lambda + M^{\, 2})^{\frac{3}{4}  }}  ,  \\   
&\textbf{(ii)} \quad  \left\|   \pi \left( \frac{ \mbf{k}^2 \tau_{\, \mbf{x}_j} \rhomLambdaa }{ \omegam^{3/2}} \right) ( \hat{\mbf{p}}_j^2  \tens \1 \; + \; M^{\, 2}  + \lambda )^{-1} (H_{0 , \, m } (\Lambda_  ) +1)^{-1} \right\| 
 \;  \leq \; 
  \gamma_j  \;  \frac{\Lambda^{2\alpha}}{\lambda + M^{\,2}} ,  \\
&\textbf{(iii)} \; \sum_{\nu =1}^3 \left\|  \pi \left( \frac{ k^\nu \tau_{\, \mbf{x}_j} \rhomLambdaa }{ \omegam^{3/2}} \right)^2( \hat{\mbf{p}}_j^2  \tens \1 \; + \; M^2  + \lambda )^{-1}  (H_{0 , \, m } (\Lambda ) +1)^{-1} \right\| 
 \;  \leq \;   \mu_j  \; \frac{\Lambda^{2\alpha}}{\lambda + M^{\, 2}} ,  
\end{align*}
respectively.
\end{lemma}
\textbf{{\large (Proof) }}  \textbf{(i)} Let $ \Lambda > 1$. It holds that for all $\Psi \in \ms{D}(H_{0 , \, m}(\Lambda ))$, 
\begin{align}
&\left\|   \pi \left( \frac{ k^\nu \tau_{\, \mbf{x}_j} \rhomLambdaa }{ \omegam^{3/2}} \right)    ( \hat{p}_j^\nu \tens \1 ) 
 ( \hat{\mbf{p}}_j^2  \tens \1 \, + \, M^{\, 2}  + \lambda )^{-1} \Psi \right\|   \notag \\
&\leq 
\left\|   \pi \left( \frac{ k^\nu \tau_{\, \mbf{x}_j} \rhomLambdaa }{ \omegam^{3/2}} \right) ( \1 \tens \Hbm +1)^{-1/2} \right\|
  \left\|  ( \hat{p}_{j}^\nu \tens \1 )  ( \hat{\mbf{p}}_j^2  \tens \1 \, + \, M^{\, 2}  + \lambda )^{-\frac{1}{4}}   (\Hpar \tens \1 +1 )^{-1/2} \right\|      \notag \\
 & \qquad \qquad \times \|  ( \hat{\mbf{p}}_j^2  \tens \1 \, + \, M^{\, 2}  + \lambda )^{-\frac{3}{4}} \| \,
 \left\| ( \1  \tens \Hbm +1)^{1/2}  (\Hpar \tens \1 +1 )^{1/2} \Psi \right\|   \notag \\
 & \leq \frac{1}{(\lambda + M^{\, 2} )^{\frac{3}{4} }}\left\|  
 \pi \left( \frac{ k^\nu \tau_{\, \mbf{x}_j} \rhomLambdaa }{ \omegam^{3/2}} \right) ( \1 \tens \Hbm +1)^{-1/2} \right\|  \left\| ( \1 \tens \Hbm +1)^{1/2}  (\Hpar \tens \1+1 )^{1/2} \Psi \right\| . \label{2014-5/9-4}
\end{align}
Note that $ \| \tau_{\, \mbf{x}} f\| = \|f \|$, $f \in \LtwoRthree $. We see that   $ \left\|   \frac{ |\mbf{k}| \rhomLambdaa }{ \omegam^{3/2}}  \right\| 
 \, \in O (\Lambda^\alpha ) $ and $  \left\|   \frac{ |\mbf{k}| \rhomLambdaa }{ \omegam^{2}}  \right\|   \, \in O (\Lambda^{ \frac{\alpha}{2}} ) $ for all $m \geq 0$. Then,  we have 
\begin{equation}
\left\| \frac{}{}  \pi \left( \frac{k^\nu \tau_{\, \mbf{x}_j} \rhomLambdaa }{ \omegam^{3/2}} \right) ( \1  \tens \Hbm +1)^{-1/2} \right\| 
\; \; \in \; O (\Lambda^\alpha ) .  \label{order1}
\end{equation}
Since  $\Lambda \geq1$, we have $ \left\| ( \1  \tens \Hbm +1)^{1/2}  (\Hpar \tens \1 +1 )^{1/2} \Psi \right\|  \leq
\left\|  H_{0 , \, m} (\Lambda ) \Psi \right\|   + \| \Psi \|  $.  By (\ref{2014-5/9-4}) and (\ref{order1}) it holds  that 
\begin{equation}
\left\|  \pi \left( \frac{ k^{\nu} \tau_{\, \mbf{x}_j} \rhomLambdaa }{ \omegam^{3/2}} \right)
  ( \hat{\mbf{p}}_j^\nu \tens \1 ) 
 (  \hat{\mbf{p}}_j^2 \tens \1 \, + \,  M^{\, 2}  + \lambda )^{-1} \Psi \right\| 
 \;  \leq \; 
\left( \beta_j \frac{\Lambda^\alpha}{(\lambda + M^{\, 2} )^{\frac{3}{4} }} \right) \left( \left\|  H_{0 , \, m } (\Lambda ) \Psi \right\|   + \| \Psi \| \frac{}{} \right) . \notag
\end{equation}
Thus  \textbf{(i)} is obtained. \\
\textbf{(ii)}  It is  seen that
\begin{align}
& \left\|  \pi \left( \frac{ \mbf{k}^2 \tau_{\, \mbf{x}_j} \rhomLambdaa }{ \omegam^{3/2}} \right) (  \hat{\mbf{p}}_j^2 \tens \1 \, + \, M^{\, 2}  + \lambda )^{-1}  \right\| \notag \\
 & \qquad \qquad   \leq
\frac{1}{M^{\, 2} + \lambda} \left\|  \pi \left( \frac{ \mbf{k}^2 \tau_{\, \mbf{x}_j} \rhomLambdaa }{ \omegam^{3/2}} \right)  
 (\1 \tens \Hb +1 )^{-1/2}  \right\| \;  \left\|  (\1 \tens \Hbm +1 )^{1/2} \Psi \right\| . \label{2014-5/9-3}
 \end{align}
Since $ \left\|   \frac{ \mbf{k}^2 \rhomLambdaa }{ \omegam^{3/2}}  \right\| 
 \; \in O (\Lambda^{2 \alpha }) $ and $  \left\|   \frac{ \mbf{k}^2 \rhomLambdaa }{ \omegam^{2}}  \right\| 
\; \in O (\Lambda^{ \frac{3}{2}\alpha } ) $ for all $m \geq 0$, it follows  that
\begin{equation}
\left\|  \pi \left( \frac{ \mbf{k}^2  \tau_{\, \mbf{x}_j}  \rhomLambdaa }{ \omegam^{3/2}} \right)    ( \1 \tens \Hbm +1 )^{-1/2}  \right\| \, \in \, O (\Lambda^{2 \alpha}). \label{order2}
\end{equation}
We see that $   \| ( \1 \tens \Hbm +1 )^{1/2} \Psi \|  \leq \| (H_{0 , \, m} (\Lambda ) + 1 ) \Psi \| $,  $\Lambda \geq 1$. Then by  (\ref{2014-5/9-3}) and (\ref{order2}),  
\[
\left\|  \pi \left( \frac{ \mbf{k}^2  \tau_{\, \mbf{x}_j}  \rhomLambdaa }{ \omegam^{3/2}} \right) ( \hat{\mbf{p}}_j^2  \tens \1 \; + \; M^{\, 2}  + \lambda )^{-1} \Psi \right\| 
 \; \leq \; 
  \gamma_j  \;  \frac{\Lambda^{2\alpha}}{\lambda + M^{\, 2}}  \left( \left\|  
H_{0 , \, m} (\Lambda ) \Psi \right\|   + \| \Psi \| \right) . 
\]
Thus  \textbf{(ii)} is proven. \\
\textbf{(iii)} It is seen that {\small
\begin{align}
 & \left\|  \pi \left( \frac{ k^{\nu}  \tau_{\, \mbf{x}_j}  \rhomLambdaa }{ \omegam^{3/2}} \right)^2 (  \hat{\mbf{p}}_j^2 \tens \1 \, + \, M^{\, 2}  + \lambda )^{-1} \Psi  \right\| 
\notag \\
& \qquad \qquad  \qquad \quad   \leq
\frac{1}{M^{\, 2} + \lambda} \left\|  \pi \left( \frac{ k^\nu  \tau_{\, \mbf{x}_j}  \rhomLambdaa }{ \omegam^{3/2}} \right)^2  
 ( \1  \tens \Hbm +1 )^{-1}  \right\| \;  \left\| (\1  \tens \Hbm +1 ) \Psi \frac{}{} \right\|
\label{2014-5/9-5} 
 \end{align}}
 From the subsequent Lemma  \ref{5/15-a}, we have{\small
 \begin{equation}
  \left\|  \pi \left( \frac{ k^\nu  \tau_{\, \mbf{x}_j}  \rhomLambdaa }{ \omegam^{3/2}} \right)^2  
 ( \1 \tens \Hbm +1 )^{-1}  \right\|    \leq    
 2  \left\|   \frac{ |\mbf{k}|  \rhomLambdaa }{ \omegam^{2}}  \right\|^2  \, 
+ \,  2   \left\|   \frac{ |\mbf{k}| \rhomLambdaa }{ \omegam^{2}}  \right\|  \times   \left\|   \frac{ |\mbf{k}| \rhomLambdaa }{ \omegam^{3/2}}  \right\|   
 \, +  \, \frac{3}{2} \left\|   \frac{ |\mbf{k}| \rhomLambdaa }{ \omegam^{3/2}}  \right\|^2  
   \notag
 \end{equation}}
Note that $ \;  \left\|   \frac{ |\mbf{k}|  \rhomLambdaa }{ \omegam^{2}}  \right\|^2 
\; \in O (\Lambda^{ \alpha } ) $,  
 $ \; \left\|   \frac{ |\mbf{k}| \rhomLambdaa }{ \omegam^{2}}  \right\|  \times \left\|   \frac{ |\mbf{k}| \rhomLambdaa }{ \omegam^{3/2}}  \right\|  \; \in O (\Lambda^{ \frac{3}{2}\alpha }  ) \;   $ and $ \left\|   \frac{ |\mbf{k}| \rhomLambdaa }{ \omegam^{3/2}}  \right\|^2   \; \in O (\Lambda^{2 \alpha }) $,
 for all $m \geq 0$. Then, we have
\begin{equation}
\left\|  \pi \left( \frac{ k^\nu \tau_{\, \mbf{x}_j}  \rhomLambdaa }{ \omegam^{3/2}} \right)^2  
 ( \1 \tens \Hbm +1 )^{-1}  \right\| \; \in \; O (\Lambda^{2 \alpha}). \label{order3}
\end{equation}
Since $    \| (\1 \tens \Hbm +1 )^{1/2} \Psi \|  \leq \| (H_{0 , \, m} (\Lambda) + 1 ) \Psi \| $,  $\Lambda \geq 1$, it follows 
 from (\ref{2014-5/9-5}) and (\ref{order3}) that
\[
\left\|  \pi \left( \frac{ k^\nu  \tau_{\, \mbf{x}_j}  \rhomLambdaa}{ \omegam^{3/2}} \right)^2( \hat{\mbf{p}}_j^2  \tens \1 \, + \, M^{\, 2}  + \lambda )^{-1}  \Psi  \right\| 
 \; \; \leq \; \;  \mu_j  \; \frac{\Lambda^{2\alpha}}{\lambda + M^{\, 2}} \left( \left\|  H_{0 , \, m} (\Lambda ) \Psi \right\|   + \| \Psi \| \right) .
\]
Thus we obtain \textbf{(iii)}. $\blacksquare $ \\

  \begin{lemma} \label{5/15-a} \upshape
  Let $\omega$ be a non-negative function on $\Rd $, $d \geq 1$, and   $f , \; g \;  \in \ms{D} (\omega^{-1/2})  $. Then for all $\Psi \in  \ms{D}(\sqzb{\omega})  $,
\begin{align*}
&\mathbf{(i)} \quad 
 \| a (f) a(g)  \Psi \| \leq \left\| \frac{f}{\sqrt{\omega}} \right\| 
 \left\| \frac{g}{\sqrt{\omega}} \right\|  \, \| \sqzb{\omega} \Psi \|   , \\
&\mathbf{(ii)} \quad
\| a^{\ast} (f) a(g)  \Psi \| \leq \left\| \frac{f}{\sqrt{\omega}} \right\|  \left\| \frac{g}{\sqrt{\omega}} \right\|  \, 
 \| \sqzb{\omega}  \Psi \| +  \| f \| \left\| \frac{g}{\sqrt{\omega}} \right\|  \, \| \sqzb{\omega}^{1/2} \Psi \|   ,  \\
& \mathbf{(iii)} \quad
\| a (f) a^{\ast}(g)  \Psi \| \leq \left\| \frac{f}{\sqrt{\omega}} \right\|  \left\| \frac{g}{\sqrt{\omega}} \right\| 
 \, \|  \sqzb{\omega}  \Psi \| +  \left\| \frac{f}{\sqrt{\omega}} \right\| 
\| g \|   \, \| \sqzb{\omega}^{1/2} \Psi \|  
 +  \| f\| \, \| g \| \;  \| \Psi \| , \\
 & \mathbf{(iv)} \quad
\| a^{\ast} (f) a^{\ast}(g)  \Psi \| \leq \left\| \frac{f}{\sqrt{\omega}} \right\|  \left\| \frac{g}{\sqrt{\omega}} \right\|  \, \|\sqzb{\omega}   \Psi \|  \notag \\
&  \qquad  \qquad \qquad \qquad \qquad \qquad  + 
 \left(  \| f \| \left\| \frac{g}{\sqrt{\omega}} \right\| +  
\left\| \frac{f}{\sqrt{\omega}} \right\| 
\| g \|    \right) 
 \, \|  \sqzb{\omega}^{1/2} \Psi \|   
 +2 \| f \| \| g\| \| \Psi \| .
\end{align*}
  \end{lemma}
\textbf{(Proof)}  \textbf{(i)} is  directly proven by the Schwartz inequality, and   \textbf{(ii)} - \textbf{(iv)} can be proven  by using \textbf{(i)},    the commutation relations  (\ref{CCR}), and the basic norm inequalities (\ref{bounda}) and (\ref{boundad}). $\blacksquare$

\begin{proposition} \label{9/13.a} \upshape $\;$ \\ 
 \textbf{(i)} Assume    $m>0$. Let   $\Lambda_{\divideontimes} >0 $. Then it holds that  
for all  $\Lambda > \Lambda_{\divideontimes}$,
 \begin{equation}
  \| \VeffkmLambda  \psi  \| \; \leq \;
c (  \kappa , N ) \left( 1 +   \frac{ \Lambda^{2 \alpha }_{\divideontimes} }{\Lambda_{\divideontimes}^{2 \alpha } -m^2} \right)  \|   \Hpar \psi \| \;  ,  \qquad \psi \in \ms{D} (\Hpar ) ,  \notag
 \end{equation}
where $ c (  \kappa , N ) = \frac{\kappa^2 }{16 \pi } N(N-1) $. In addition,  
\begin{equation}
\text{s}-\lim_{\Lambda \to \infty }  \VeffkmLambda \psi 
 \;  = \;    \Veffkm  \psi , \qquad \psi \in \ms{D} (\Hpar )  . \notag
\end{equation}
$\;$ \\
\textbf{(ii)} Let $\VeffkLambda = \VeffkzeroLambda$ and $\Veffk = \Veffkzero $. Then 
for all  $\Lambda > 0$,
 \begin{equation}
  \qquad \| \VeffkLambda  \psi  \| \; \leq \;
2 \, c (\, \kappa , N )  \|   \Hpar \psi \| \;  ,  \qquad  \qquad \psi \in \ms{D} (\Hpar ) . 
\notag
 \end{equation}
 Moreover, 
\begin{equation}
\qquad \text{s}-\lim_{\Lambda \to \infty }  \VeffkLambda  \psi
 \;  = \;  \Veffk   \psi , \qquad \psi \in \ms{D} (\Hpar )  . \qquad \qquad  \notag
\end{equation}
\end{proposition}

\begin{remark} \upshape
For  Proposition \ref{9/13.a} \textbf{(i)},
we see that $c (\kappa , N ) < \frac{1}{2} $, $ 0 <  |\kappa| < 
 \sqrt{\frac{8  \pi}{ N (N-1) }}  $. Then  we can take $\epsilon_{\divideontimes} > 0$ such as $ 2 + \epsilon_{\divideontimes} < \frac{1}{c (\kappa , N )} $. Since 
$\limsup\limits_{\Lambda \to \infty} \frac{ \Lambda^{2 \alpha } }{\Lambda^{2 \alpha } -m^2}  =1$, there exists $\Lambda (\epsilon_{\divideontimes}) > 0$ such that for all 
$\Lambda \geq  \Lambda ( \epsilon_{\divideontimes} )$, it holds that  $   \frac{ \Lambda^{2 \alpha } }{\Lambda^{2 \alpha } -m^2}  < 1 + \epsilon_{\divideontimes}$.
Let $ \tilde{\Lambda}_{\divideontimes}  = \max \left\{ \Lambda_{\divideontimes} , \Lambda ( \epsilon_{\divideontimes} ) \right\}$. 
Then it holds that 
$c (\kappa , N) \left( 1 +   \frac{ \tilde{\Lambda}_{\divideontimes}^{2 \alpha } }{\tilde{\Lambda}_{\divideontimes}^{2 \alpha } -m^2} \right) < 1$.  
\end{remark}
\textbf{(Proof)} 
\textbf{(i)}  Let $\Lambda > \Lambda_{\divideontimes}$. From the spherical symmetry of  $\rho_{m, \Lambda} 
(\mbf{k}) = \chi_{{}_{J_{m, \Lambda}}} (|\mbf{k}|)$,  it follows that 
  $ \int_{\Rthree} \frac{\rhomLambdaa (\mbf{k})}{\omegam( \mbf{k})^2} e^{\, i\mbf{k \cdot x} } d \mbf{k}  \; = \;   \frac{ \pi}{|\mbf{x}|} \int_{-\Lambdaa}^{\Lambdaa}
 \frac{r}{r^2 +m^2}  \sin (|x| r) d r$.
Then, 
\begin{equation} 
\VeffkmLambda
 \; \; 
 = \; \;  - \frac{\kappa^2}{16 \pi^2} 
   \sum_{j < l}^N   \; \frac{1}{|\mbf{x}_j - \mbf{x}_l |} 
 G_{m,\, \Lambdaa} ( \mbf{x}_j - \mbf{x}_l ) ,   \label{5/7.6}
\end{equation}
where $ 
 G_{ m , \, \Lambdaa } (\mbf{x} ) \; 
=  \; \text{Im } \, \int_{-\Lambdaa}^{\Lambdaa}
 \frac{r}{r^2 +m^2} e^{\, i |\mbf{x} | r} dr $. 
  Let   
  $\Gamma_{\Lambdaa}=  \{ z=t \left| \frac{}{} \right. -\Lambdaa \leq t \leq \Lambdaa\} \cup C_{\Lambdaa}$  with 
$ C_{\Lambdaa} =\{ z=\Lambdaa e^{\, i \theta} \left| \frac{}{} 0 \right. \leq \theta < \pi  \}$.
By the residue theorem,   it holds that 
  $  \int_{\Gamma_{\Lambdaa}} 
 \frac{z}{z^2 +m^2} e^{\, i |\mbf{x}| z } d z = i \pi e^{-m | \mbf{x} |} $.
Then, 
\begin{equation}
 G_{ m , \, \Lambdaa } ( \mbf{x}  ) \; = \;
 \pi e^{- m|\mbf{x}  |} \, + \, \text{Im} \int_{C_{\Lambdaa}}   \frac{z}{z^2 +m^2} 
 e^{\, i  |\mbf{x}  | z } d z .  \label{Feb5.3}
\end{equation}
We see that  
\[
 \left|   \int_{C_{\Lambdaa}}   \frac{z}{z^2 +m^2} e^{\, i |\mbf{x}| z } d z \right|
 \; \leq \;  \frac{\Lambda^{2 \alpha}}{\Lambda^{2\alpha} -m^2} 
  \int_{0}^\pi e^{-|\mbf{x} | \Lambda^{\alpha} \sin \theta } d \theta 
  \; \leq \;  \pi \frac{\Lambda^{2 \alpha}}{\Lambda^{2\alpha} -m^2} .
\]
 From (\ref{Feb5.3}), we have
\begin{equation}
\left| G_{ m ,\, \Lambdaa } (\mbf{x}_j - \mbf{x}_l )   \; \right|
 \;   
  \leq  \; \pi \left(1 +  \frac{\Lambda^{2 \alpha}}{\Lambda^{2\alpha} -m^2} \right) 
  \leq  \; \pi \left(1 +  \frac{\Lambda_{\divideontimes}^{2 \alpha}}{\Lambda^{2\alpha}_{\divideontimes} -m^2} \right)
\label{5/8.4}
\end{equation}
By (\ref{5/8.4}) and $\| \frac{1}{| \mbf{x}_j - \mbf{x}_{l} |} \psi \|_{\LtwoRthreeNx} \leq 2 \| \nabla_{\mbf{x}_{i}} \psi \|_{\LtwoRthreeNx} $, it holds that  for all $ \psi \in C_{0}^{\, \infty} (\RthreeN )$,
{\small
\begin{equation}
\left\| \VeffkmLambda \psi \right\|  \; \leq \; 
\frac{\kappa^2 }{8 \pi }  \sum_{j<l}  \left\| \frac{1}{|\mbf{x}_j - \mbf{x}_l |^2}
 G_{ m, \,  \Lambdaa } (\mbf{x}_j - \mbf{x}_l ) \psi \right\|  \leq \;  \frac{\kappa^2  N(N-1)}{16 \pi } \left( 1 +   \frac{ \Lambda_{\divideontimes}^{2 \alpha } }{\Lambda_{\divideontimes}^{2 \alpha } -m^2} \right) \|  \Hpar \psi\|.
\label{11/29.1}
\end{equation}
}Note that   (\ref{11/29.1}) follows for all  $\psi \in \ms{D} (\Hpar)$, since $ C_{0}^{\, \infty} (\RthreeN )$ is a core of $\Hpar $.
By the bounded convergence theorem,  it follows that  for all $ \mbf{x} \ne \mbf{0}$, 
\begin{equation}
\qquad \qquad  \lim_{\Lambda \to \infty} \left|   \int_{C_{\Lambdaa}}   \frac{z}{z^2 +m^2} e^{i |\mbf{x}| z } d z \right|
 \leq  \lim_{\Lambda \to \infty}   
  \int_{0}^{\pi} e^{-|\mbf{x} | \Lambda^{\alpha} \sin \theta } d \theta =0   , \notag
\end{equation}
 and we have
\begin{equation}
  \qquad \qquad
\lim_{\Lambda \to \infty  } \left| \frac{}{} G_{ m, \, \Lambdaa } (  \mbf{x}_j - \mbf{x}_l    ) \;  -   \; \pi e^{-m | \mbf{x}_j - \mbf{x}_l |}   \; \right| \; \; = \; \; 0 ,  \quad \quad \text{a.e. } (\mbf{x}_1 , \cdots, \mbf{x}_N ) \, \in \, \RthreeN . 
 \notag
\end{equation}
Note that  
   $\| \frac{1}{|\mbf{x}_i - \mbf{x}_l |} \psi  \|_{L^2} < \infty $  for all 
$\psi \in C_{0 }^{\, \infty} (\RthreeN)$. Then  the Lebesgue dominated convergence  theorem  yields that for all $\psi \in C_{0 }^{\, \infty} (\RthreeN)$, 
\begin{equation}
\lim_{\Lambda \to \infty}  \int_{\RthreeN} 
\left|  \frac{\pi e^{-m |\mbf{x}_i - \mbf{x}_l | }  - \, G_{m , \,\Lambdaa} (\mbf{x}_i - \mbf{x}_l )}{|\mbf{x}_j  - \mbf{x}_l |} \right|^2  
   \left| \frac{}{} \psi (\mbf{x}_1 , \cdots , \mbf{x}_N ) \right|^2 
 d \mbf{x}_1  \cdots d \mbf{x}_N  \; \; = \; \; 0  . \notag
 \end{equation}
Thus, we obtain for all $\psi \in C_{0 }^{\, \infty} (\RthreeN) $
\begin{equation}
\text{s}-\lim_{\Lambda \to \infty } \VeffkmLambda \psi \; \; = \; \;
  \Veffkm \psi . \label{5/9.1}
 \end{equation} 
 Since $C_{0 }^{\, \infty} (\RthreeN)$ is a core of $\Hpar $, it is  proven that  (\ref{5/9.1}) holds for all $\psi \in \ms{D} (\Hpar )$.   \\
\textbf{(ii)}  
Let $\rho_{ \Lambda} (\mbf{k}) = \rho_{0, \Lambda} 
(\mbf{k}) = \chi_{{}_{J_{0, \Lambda}}} (|\mbf{k}|)$. 
 It holds that 
  $ \int_{\Rthree} \frac{\rhoLambdaa (\mbf{k})}{|\mbf{k}|^2} e^{i\mbf{k \cdot x} } d \mbf{k}  \; = \;  
 \frac{ \pi}{|\mbf{x}|} \int_{I_{\Lambdaa }}
 \frac{1}{r}  \sin (|x| r) d r$ where $ I_{\Lambda^{\alpha}} = \{ x=t \left| \frac{}{} \right. -\Lambdaa \leq t \leq -\Lambda^{-\alpha} , \, \Lambda^{-\alpha} \leq t \leq \Lambdaa\}$.
Then we have 
\begin{equation} 
\VeffkLambda
 \; \; 
 = \; \;  - \frac{\kappa^2}{16 \pi^2} 
   \sum_{j < l}^N   \; \frac{1}{|\mbf{x}_j - \mbf{x}_l |} 
 G_{\Lambdaa} ( \mbf{x}_j - \mbf{x}_l ) ,   \label{5/7.6a}
\end{equation}
where 
$  G_{\Lambdaa } (\mbf{x}) \, 
= \,  \text{Im } \, \int_{I_{\Lambda^{\alpha}}}
 \frac{1}{r} e^{\, i |\mbf{x} | r} dr $. 
  Let   
  $\tilde{\Gamma}_{\Lambdaa}=  I_{\Lambda^{\alpha}} \cup C_{\Lambdaa} \cup \tilde{C}_{\Lambda^{-\alpha}}$  where   
$ C_{\Lambdaa} =\{ z=\Lambdaa e^{\, i \theta} \left| \frac{}{} 0 \right. \leq \theta < \pi  \}$ and 
$ \tilde{C}_{\Lambdaa} =\{ z= \frac{1}{\Lambda^{\alpha}} e^{-i \theta} \left| \frac{}{}  \right. \pi \leq \theta < 2 \pi  \}$. 
Cauchy's integral formula yields that  $  \int_{\tilde\Gamma_{\Lambdaa}} 
 \frac{1}{z} e^{\, i |\mbf{x}| z } d z = 0 $.
Then 
\begin{equation} 
 G_{\Lambdaa } ( |\mbf{x}_j - \mbf{x}_l |) 
= - \text{Im} \left\{ \int_{ C_{\Lambda^{\alpha}}} 
 \frac{1}{z} e^{\, i |\mbf{x}_j - \mbf{x}_l |  z } d z +  \int_{\tilde{C}_{\Lambda^{-\alpha}}} 
 \frac{1}{z} e^{\, i |\mbf{x}_j - \mbf{x}_l | z } d z     \right\}    \label{Feb5.7}
\end{equation}
It is seen that 
\[
 \left|   \int_{C_{\Lambdaa}}   \frac{1}{z} e^{\, i |\mbf{x}| z } d z \right|   
 \; \leq \;  \int_{0}^{\pi} e^{- |\mbf{x}|\Lambda^{\alpha} \sin \theta } d \theta \; \leq \; \pi  , \qquad \left|   \int_{\tilde{C}_{{\Lambda}^{-\alpha}}}   \frac{1}{z} e^{i |\mbf{x}| z } d z \right|   
\; \leq \;  \int_{\pi}^{2 \pi} e^{ |\mbf{x}|\Lambda^{-\alpha} \sin \theta } d \theta \; \leq  \;  \pi  . 
\]
Then we have $\left| G_{\Lambdaa } ( \mbf{x}_j - \mbf{x}_l ) \right| \leq 2 \pi  $.  By  (\ref{Feb5.7}) and $\| \frac{1}{| \mbf{x}_j - \mbf{x}_{l} |} \psi \|_{\LtwoRthreeNx} \leq 2 \| \nabla_{\mbf{x}_{i}} \psi \|_{\LtwoRthreeNx} $, it follows that   for all  $ \psi \in C_{0}^{\, \infty} (\RthreeN )$,
\begin{equation}
\left\| \VeffkLambda \psi \right\|  \; \leq \; \frac{\kappa^2}{8 \pi } N (N-1) \|  \Hpar \psi\|. \label{Feb6.1}
\end{equation}
 Note that  (\ref{Feb6.1}) follows for all  $\psi \in \ms{D} (\Hpar)$ since $ C_{0}^{\infty} (\RthreeN )$ is a core of $\Hpar  $. 
By (\ref{Feb5.7}), we have
\begin{equation} 
\left|  G_{\Lambdaa } ( \mbf{x}_j - \mbf{x}_l ) - \pi \, \right|
 \; \leq \;   \left| \int_{ C_{\Lambda^{\alpha}}} 
 \frac{1}{z} e^{\, i |\mbf{x}_j - \mbf{x}_l |  z } d z \right| \,   + \, \left| \int_{\tilde{C}_{\Lambda^{-\alpha}}} 
 \frac{1}{z} e^{\, i |\mbf{x}_j - \mbf{x}_l | z } d z   \, + i \pi \, \right|  .  \notag
\end{equation} 
By the bounded convergence theorem, it follows that  for all $\mbf{x} \in \mbf{\Rthree} \backslash \{ \mbf{0} \} $, 
\begin{align}
& \lim_{\Lambda \to \infty} \left| \int_{ C_{\Lambda^{\alpha}}} 
 \frac{1}{z} e^{\, i |\mbf{x} |  z } d z \right|  \leq 
  \lim_{\Lambda \to \infty}   \int_{0}^{ \pi }  e^{ -|\mbf{x}| \Lambda^{\alpha } \sin \theta   }   d \theta   \; = 
\; 0 , \notag  \\
& \lim_{\Lambda \to \infty} \left|
\int_{\tilde{C}_{\Lambda^{-\alpha}}}  \frac{1}{z} e^{i |\mbf{x}| z } d z \, + \, i \pi \, 
 \right|  \leq  \lim_{\Lambda \to \infty}   \int_{\pi}^{2 \pi } 
\left| e^{ \, |\mbf{x}| \Lambda^{-\alpha } (\sin \theta + i \cos \theta)  } -1  \right| d \theta  \; = \; 0  . \notag
\end{align}
Then  we have
\begin{equation}
  \quad \quad
\lim_{\Lambda \to \infty  } \left|  G_{\Lambdaa } (\mbf{x}_j - \mbf{x}_l )  \, -  \, \pi   \; \right| 
\;   = \;  0 ,  \quad \quad \text{a.e. } (\mbf{x}_1 , \cdots, \mbf{x}_N ) \, \in \, \RthreeN . 
\notag 
\end{equation}
By the Lebesgue dominated convergence  theorem, we have   for all $\psi \in C_{0 }^{\, \infty} (\RthreeN)$, 
\begin{equation}
\lim_{\Lambda \to \infty}  \int_{\RthreeN} 
\left| \frac{  \pi - \, G_{\Lambdaa} (\mbf{x}_i - \mbf{x}_l )}{|\mbf{x}_j  - \mbf{x}_l |} \right|^2  
   \left|  \psi (\mbf{x}_1 , \cdots , \mbf{x}_N ) \right|^2 
 d \mbf{x}_1  \cdots d \mbf{x}_N  \;  =  \; 0  .  \notag 
 \end{equation}
From this,   it is  proven  that 
$\text{s}-\lim\limits_{\Lambda \to \infty } \VeffkLambda \psi \, = \,
  \Veffk \psi $ for all $\Psi \in \ms{D} (\Hpar ) $, since $C_{0 }^{\, \infty} (\RthreeN)$ is a core of $\Hpar $.  $\blacksquare$

\begin{lemma} \label{u-convergence} \upshape  Let $0< \alpha < 2$. 
Then 
\begin{equation}
\text{s}-  \lim_{\Lambda \to \infty}  \UkmLambda  \; = \; \1 . \notag 
\end{equation}
\end{lemma}
(\textbf{Proof})
Let $\Psi \in \corem $. 
 We see that 
\begin{equation}
\left\|  
 \pi \left( \frac{ \tau_{\, \mbf{x}}  \rhomLambdaa}{\sqrt{\omegam }^{\, 3}} \right) \Psi \right\|
 \leq \left\| \frac{\rhomLambdaa}{\omegam^{\,2} } \right\| \,  
\left\| \left( \1 \tens \Hbm^{1/2} \right) \Phi  \right\|  \;  + \;
\left\|  \frac{\rhomLambdaa}{\sqrt{\omegam }^{\, 3}}  \right\| \, \left\| \Psi \right\|  .
\notag \end{equation}  
Note $ \left\|  \frac{\rhomLambdaa}{\omegam^{\, 2} }   \right\|  \in O (\Lambda^{{\alpha} /2})$ and
$ \left\|  \frac{\rhomLambdaa}{\sqrt{\omegam }^3}  \right\|  \in O ( \log \Lambda ) $ for all $m \geq 0 $. We see that 
$|e^{\,i \lambda }-1| \leq |\lambda | $, $\lambda \in \mbf{R}$. Then it holds that 
{\small
\begin{equation}
\left\|  \left( \UkmLambda - \1 \right) \Psi \right\| 
 \; \leq \; \frac{|\kappa|}{\Lambda} \sum_{j=1}^N \left\|  
 \pi \left( \frac{ \tau_{\, \mbf{x}_j}  \rhomLambdaa}{\sqrt{\omegam }^{\, 3}} \right) \Psi \right\| \;  
\leq \;  c |\kappa| N \left(  \frac{1}{\Lambda^{1-\frac{\alpha}{2}}} \left\| \left( \1 \tens \Hbm^{1/2} \right) \Psi  \right\|
 +  \frac{\log \Lambda }{\Lambda } \left\| \Psi \right\| \right)  , \notag
\end{equation}
}where $c \geq 0$ is a constant. Note that  $ \lim\limits_{\Lambda \to \infty } \Lambda^{\frac{\alpha}{2}-1}=0 $,  $0 <\alpha< 2$, and $ \lim\limits_{\Lambda \to \infty } \frac{\log \Lambda }{\Lambda }=0 $. Then we have s$- \lim\limits_{\Lambda \to \infty} \UkmLambda \Psi = \Psi $.
Since $\corem$ is  dense in $\ms{H}$ and $\|\UkmLambda \| =1$,  the assertion is proven. $\blacksquare $

$\quad $ \\
{\Large \textbf{(Proof of Theorem \ref{MainTheorem})}} \\ 
By the dressing transformation, it holds that       for all $ z \in \mbb{C} \backslash \mbb{R} $, 
\begin{equation}
 \UkmLambda^{-1}  \left(  \HkmLambdaren   -  z \right)^{-1} 
 \UkmLambda 
 \, = \, \left(
  H_{0, m } (\Lambda )   
    +  \VeffkmLambda   +  \deltaHkmLambda    -  z  \right)^{-1}. 
 \end{equation}
  By Lemma \ref{u-convergence},  we have s-$\lim\limits_{\Lambda \to \infty} \UkmLambda = \1$. From Lemma \ref{A.3'}, Proposition \ref{5/7.5} and Proposition  \ref{9/13.a}, we  see that  $  H_{0, m } (\Lambda )  \, 
    + \,  \VeffkmLambda \,  + \, \deltaHkmLambda    $ satisfies 
\textbf{[A.1]} - \textbf{[A.3]} with applying $H_0 (\Lambda ) $ to $X_0 (\Lambda^2 )$, 
  $ \deltaHkmLambda $ to $C_{\I} (\Lambda^2 )$, $0$ to $C_{\I}$, 
$\VeffkmLambda $ to $C_{\II} (\Lambda^2 ) $ and $\Veffkm $ to $C_{\II}$. 
 Then we have
\begin{align*}
  \text{s} \, -\, \lim_{\Lambda \to \infty }\left(  \HkmLambdaren   -  z \right)^{-1} 
  \, &=   \text{s}  -  \lim_{\Lambda \to \infty } \left( 
  H_{0, m } (\Lambda )       +   \VeffkmLambda   +  \deltaHkmLambda    -  z  \right)^{-1} \\
&  = \left( \Hpar \tens \1 + \Veffkm \tens \Pb -z \right)^{-1} (\1 \tens \Pb ) \\
& =  \left( \Hpar  + \Veffkm  -z \right)^{-1}  \tens \Pb .
\end{align*}  
Thus proof is finished. $\blacksquare $

$\quad$ \\
{\Large Acknowledgments} $\; $  This research is supported by JSPS KAKENHI Grant Numbers JP16K17607. 

\end{document}